%% file: main.tex
\def\href#1{}
\documentclass[12pt]{amsart}

\input{PlainPaper/Package/script}
\input{PlainPaper/Command/script}

\begin{document}
\input{PlainPaper/Style/script}
\author{Qirui Li}
\title{On Gross-Keating's result of lifting endomorphisms for formal modules}
\date{}
\maketitle
\input{Lifting/Command/script}

\def\load#1{\input{Lifting/#1}}
\begin{abstract}
\input{PlainPaper/Command/Abstract}

\end{abstract}
\tableofcontents
\def\use#1{}
\def\s#1{}
\def\Road#1{}
\section{Introduction}
\input{PlainPaper/Command/Introduction}
\input{PlainPaper/Command/ContentList}

\bibliography{Bib/bib}
\bibliographystyle{amsalpha}
\end{document}

%% file: PlainPaper/Package/script.tex
\usepackage{pdfsync}
\usepackage{fullpage}
\usepackage{mathrsfs}
\usepackage{latexsym, amssymb, amsmath, amscd, amsthm, amsxtra}
\usepackage{eucal}
\usepackage[utf8]{inputenc}
\usepackage{verbatim}
\usepackage[english]{babel}
\usepackage{diagbox}
\usepackage{mathtools}
\usepackage{todonotes} 
\usepackage{cite}
\usepackage{graphicx}
\usepackage{amssymb}
\usepackage{epstopdf}
\usepackage{tikz}
\usetikzlibrary{calc}
\usetikzlibrary{graphs}
\usepackage[all]{xy}
\usepackage{color}
\usepackage{fancyhdr}
\usepackage{arydshln}
\usepackage{times}
\usepackage{datetime}
\usepackage{scrtime}
\usepackage{CJKutf8}
\usepackage{cancel}
\usepackage{ulem}
\usepackage{bbm}

%% file: PlainPaper/Command/script.tex
\def\load#1{\input{PlainPaper/Command/#1}}
\input{PlainPaper/Command/basic}

\input{PlainPaper/Command/symbol}
\input{PlainPaper/Command/homologic}
\input{PlainPaper/Command/linear}

%% file: PlainPaper/Command/basic.tex
\long\def\[#1]#2{ \begin{#1}#2\end{#1}}

\renewcommand{\(}{\left(}
\renewcommand{\)}{\right)}
\renewcommand{\[}{\left[}

\newcommand{\QQ}{\mathbb{Q}}
\newcommand{\PP}{\mathbb{P}}

\newcommand{\map}[3]{#1:#2\longrightarrow #3}

\newcommand{\maps}[5]{\begin{array}{llrll}
#1&:&#2&\longrightarrow &#3\\\\
&&#4&\longmapsto &#5\\
\end{array}}

\newcommand{\dd}{\mathrm{d}}

\newcommand{\Gal}{\mathrm{Gal}}
\newcommand{\CFF}[1]{\overline{\FF_#1}}
\newcommand{\FF}{\mathbb{F}}

\newcommand{\OO}[1]{\mathcal{O}_{#1}}
\newcommand{\CO}[1]{\mathcal{O}_{\breve{#1}}}

\newcommand{\vv}{\textbf{v}}

\newcommand{\ep}{\boldsymbol{\epsilon}}

\DeclareMathOperator{\End}{\mathrm{End}}
\newcommand{\Vol}{\mathrm{Vol}}

\DeclareMathOperator{\Spf}{\mathrm{Spf}\text{ }}

\DeclareMathOperator{\length}{\mathrm{length}}
\DeclareMathOperator{\Aut}{\mathrm{Aut}}

\DeclareMathOperator{\GL}{\mathrm{GL}}


%% file: PlainPaper/Command/homologic.tex
\newcommand{\ezer}[4]{\xymatrix{
{#1}\ar@<-.5ex>[rr]_{#3}\ar@<.5ex>[rr]^{#4}&&{#2}}
}

%% file: PlainPaper/Command/linear.tex
\newcommand{\mm}[4]{\begin{pmatrix}
#1 & #2\\
#3 & #4\\
\end{pmatrix}}

\newcommand{\rr}[2]{\begin{pmatrix}
#1 & #2\\
\end{pmatrix}}

\newcommand{\cc}[2]{\begin{pmatrix}
#1\\
#2\\
\end{pmatrix}}

%% file: PlainPaper/Style/script.tex
\def\Load#1{\input{PlainPaper/Style/#1}}
\let\yourlabel=\label
\def\lb#1{\yourlabel{#1}\comments{{\color{blue}\textbf{Lable}}:\textbf{#1   }  }}
\def\status{\use{Status}}
\def\test{t}
\def\comments#1{\ifthenelse{\equal\status\test}
{{\color{red} #1}}
{}
}
\long\def\s#1{\ifthenelse{\equal\status\test}
{{\color{blue} \textbf{{\color{red}]}#1 \color{red}[}}}
{}
}

\long\def\ss#1{\ifthenelse{\equal\status\test}
{{\color{purple} \textbf{{\color{red}\}}#1 \color{red}\{~\\}}}
{}
}

\newcommand{\commenta}[1]{{\color{blue}(\textbf{Comments:} #1)}}
\newcommand{\ada}[1]{#1}
\newcommand{\dla}[1]{}
\newcommand{\ria}[1]{}
\newcommand{\xo}[1]{\input{Content/#1}}
\newcommand{\oo}[1]{\input{Content/#1}}
\newcommand{\ox}[1]{}
\input{PlainPaper/Style/plain}

%% file: PlainPaper/Style/plain.tex
\newtheorem{thm}{Theorem}[section]
\newtheorem{pof}{Proof}[section]
\newtheorem{defi}[thm]{Definition}
\newtheorem{prop}[thm]{Proposition}
\newtheorem{cor}[thm]{Corollary}
\newtheorem{summ}{Summary}
\newtheorem{rem}[thm]{Remark}
\newtheorem{ex}[thm]{Problem}
\newtheorem{lem}[thm]{Lemma}
\newtheorem{exa}[thm]{Example}
\newtheorem{proj}{Project}
\newtheorem{conj}{Conjecture}
\numberwithin{equation}{section}
\setcounter{tocdepth}{1}    

%% file: Lifting/Command/script.tex
\def\load#1{\input{Lifting/Command/#1}}
\input{PlainPaper/Command/formalgroups}

\input{PlainPaper/Command/liftingendomorphism}
\input{PlainPaper/Command/convenient}

%% file: Lifting/Command/formalgroups.tex
\newcommand{\GG}{\mathcal{G}}

%% file: Lifting/Command/liftingendomorphism.tex
\newcommand{\GK}{\GG_K}
\newcommand{\GU}{\GG^{\mathrm{univ}}}
\newcommand{\GO}{\GG_0}
\newcommand{\kw}{K_{(w)}}
\newcommand{\kwo}{K_{(w+1)}}
\newcommand{\W}{W}
\renewcommand{\O}{\mathcal{O}}
\newcommand{\e}{e}
\newcommand{\w}{w}
\newcommand{\DEF}{\mathcal{M}}
\newcommand{\blin}{\tau}
\newcommand{\biso}{\varphi}
\newcommand{\Nrd}{\mathrm{Nrd}}
\newcommand{\va}{v}
\newcommand{\y}{y}
\newcommand{\yy}{y}
\newcommand{\pp}{p}
\newcommand{\p}{p}
\newcommand{\x}{x}
\newcommand{\z}{z}
\renewcommand{\a}{{\overline{z}}}
\newcommand{\ACj}{\gamma}
\newcommand{\normal}{normal }

%% file: Lifting/Command/convenient.tex
\long\def\[#1]#2{ \begin{#1}#2\end{#1}}
\newcommand{\Proj}{\mathrm{Proj}}

%% file: Lifting/Abstract.tex
Let $K/F$ be a quadratic extension of non-Archimedean local fields of characteristic not equal to 2, with rings of integers denoted by $\OO K$ and $\OO F$. We consider a formal $\OO F$-module $\GG$, over a discrete valuation ring $\OO W$ with an uniformizer $\varpi$, with extra endomorphisms by a subring $\O$ of $\OO K$, and the height of its reduction $\GG_0=\GG\otimes \OO W/\varpi$ is 2. The endomorphism ring of $\GG_n=\GG\otimes \OO W/\varpi^{n+1}$ is a subring between $\OO s$ and $\OO D=\End(\GG_0)$. We will determine them explicitly. This result was previously proved by Gross and Keating. Their treatment is the formal cohomology theory. We will give another proof using the intersection formula of CM cycles in Lubin-Tate deformation spaces.

%% file: Lifting/Introduction.tex
\s{Goal}This article is an application of the author's intersection formula \cite{qirui2017love} to compute endormophism rings of canonical and quasi-canonical liftings. These results was firstly calculated by Gross in Proposition 3.3 of \cite{gross1986canonical} for the canonical lifting case. Then Keating treated general quasi-canonical lifting cases in Theorem 5.1 of \cite{keating1988lifting}. Their treatment is the formal cohomology theory. In contrast, we give another proof of these results using the intersection formula of CM cycles in Lubin-Tate spaces. Let $K/F$ be a separable quadratic extension, $\FF_q$ the residue field for $F$. Let $\zeta$ be an $\OO F$-algebra generator of $\OO K$ if $K/F$ is unramified, or an uniformizer of $\OO K$ if $K/F$ is ramified. Let $w$ be a non-negative integer and take $\mu=\pi^\w\zeta$. Let $\O=\OO F\oplus\mu\OO F$ be an $\OO F$-subalgebra of $\OO K$. Let $\GG_0$ be a formal $\OO F$-module of height 2 over $\CFF q$ and $\GG$ be a (quasi-)canonical lifting of $\GG_0$ as a formal $\O$-module over the abelian extension $W$ of $K$ with norm subgroup $\O^\times\subset K^\times$. Let $\varpi$ be an uniformizer of $\OO W$. Our main result is a proof of Theorem 5.1 of \cite{keating1988lifting}. Keating present his results into a list of formulae based on different cases. In contrast, we could write them uniformly as an integral. 
\begin{thm}[Equivalent to Theorem 5.1 of \cite{keating1988lifting}]\label{buxiang}\comments{buxiang}
 Suppose $q$ is odd and $\ACj$ is an automorphism of $\GG_0$ as a formal $\OO F$-module, then it can be lifted to $\GG_n$ if and only if
$$
n<\int_{\O^\times}|x-\ACj|_D^{-1}\dd x.
$$
Here $\dd x$ is the normalized Haar measure of $\O^\times$.
\end{thm}
\s{Content}Any automorphism $\ACj$ can be lifted as an isomorphism. This Theorem determines the maximal $n$ where $\ACj$ can be lifted as an automorphism. In Section \S\ref{zhunbei}, we see $\GG$ induces a closed embedding $\yy$ of $\Spf \OO W$ to the Lubin-Tate deformation space $\DEF$ as defined in \eqref{xiangxiang}. In Section \S\ref{npy}, we prove the main Theorem by our intersection formula. 

%% file: Lifting/ContentList.tex
\Road{Introduction}
\Road{Abstract}
\Road{Acknowledgement}
\Road{Command/script}
\Road{etc}
\input{PlainPaper/Command/endo}

\input{PlainPaper/Command/proofs}

\input{PlainPaper/Command/relation}

%% file: Lifting/endo.tex
\section{The endomorphism ring of formal modules}\lb{zhunbei}

\ss{This section sets up all background knowledge and notation. \[itemize]{
\item Notation: $\GG,\GK,\GO,\psi$
\item The *-isomorphism and parameter for the deformation
\item The action of division algebra
\item The action of Galois group
\item The deformation space and the meaning of intersection number
}}

In this section, we introduce canonical and quasi-canonical liftings, and relate their endomorphism rings to intersection numbers in Lubin-Tate deformation spaces.

\ss{\[itemize]{
\item What is a formal module; What is the height;
\item What is the canonical lifting; What is a quasi-canonical lifting
\item Regarding $\O$ as the subring of endomorphism ring of $\GO$.
\item What is a deformation, what is a parameter for a deformation.
\item Action of the maximal order of division algebra and action of Galois group on parameters.
}}
\subsection{Canonical and quasi-canonical liftings} We review the notion of canonical and quasi-canonical liftings follow the paper \cite{gross1986canonical} of Gross.
\ss{What is a formal module; What is the height;}

\subsubsection{Formal modules}A formal $\OO F$-module $\GG$ over a Notherian complete $\OO F$-algebra $A$ is an one dimensional formal group law over $A$ with endomorphisms by $\OO F$, such that the induced action of $\OO F$ on the Lie algebra of $\GG$ agrees with the structure map. Fix a coordinate of $\GG$ and let $[k]_\GG(X)$ be the power series over $A$ defining the action of $k\in\OO F$. If $A$ is of characteristic p, then $[\pi]_\GG(X)=\beta(X^{q^h})$ for some $\beta^\prime(X)\neq 0$, we call this $h$ the height of $\GG$.
\ss{What is the canonical lifting;What is a quasi-canonical lifting}
\subsubsection{The Canonical lifting}
For any $\GG_0$ a formal $\OO K$-module of height 1 over $\CFF q$, Lubin and Tate constructs in \cite{lubin1965formal} a formal $\OO K$ module $\GK$ over $\CO K$ such that its special fiber is $\GG_0$. They showed that this $\GK$ is unique up to isomorphism, we call $\GK$ the canonical lifting of $\GG_0$. Let $\OO D := \End(\GG_0)$ be the endomorphism ring of $\GG_0$ as a formal $\OO F$-module. Note $\GG_0$ has height 2 as a formal $\OO F$-module. Our $\End(\GG_0)$ is a maximal order of the quaternion algebra $D$ over $F$. 
\subsubsection{Quasi-canonical liftings}\lb{yuanwang}
In general, fix an embedding $\O\subset\OO D$, we can construct a lifting $\GG$ of $\GG_0$ such that the endormorphism ring $\End(\GG)\subset\End(\GG_0)$ is exactly given by $\O\subset \OO D$. The construction is made as following. The Tate module $T\GK$ of $\GK$ is a free $\OO K$-module of rank 1. Let $\vv$ be an $\OO K$ generator of $T\GK$ and the submodule $T^\prime=\O\cdot \vv\subset T\GK$ give rise to an isogeny 
\begin{equation}\lb{fangong}\map{\Phi}{\GK}{\GG}\end{equation}
of formal $\OO F$-modules with the kernel isomorphic to $\OO K/\O$ and $T\GG=T^\prime$. 
Since $\mu\cdot T\GG\subset T\GG$, $\mu$ induces an endomorphism on $\GG$. Therefore $\GG$ is a formal $\O$-module. Since all height 2 formal $\OO F$-modules over $\CFF q$ are isomorphic, we can find an isomorphism $\psi$ such that the reduction of $\psi\circ\Phi$ is an endomorphism of $\GG_0$. Replacing $\Phi$ by $\psi\circ\Phi$ we can assume $\GG$ reduces to $\GO$ over $\CFF q$. The reduction of actions of $\O$ and $\OO K$ on $\GG$ and $\GK$ induces natural embeddings
\[equation]{
\lb{zhyy}\map{\overline{[\cdot]}_{\GG}}\O{\End(\GO)}}
\[equation]{
\lb{zuohe}\map{\overline{[\cdot]}_{\GK}}{\OO K}{\End(\GO)}}
Let $\overline{[\Phi]}$ be the reduction of the map $\Phi$. For any $k\in\O$, we have
\[equation]{\lb{youyin}
\overline{[\Phi]}\circ\overline{[k]}_{\GK}=\overline{[k]}_{\GG}\circ\overline{[\Phi]}
}
In fact, the $\End(\GO)$ is the maximal order of the quaternion algebra over $F$. We can write the above equation by $\overline{[\Phi]}\circ\overline{[k]}_{\GK}\circ\overline{[\Phi]}^{-1}=\overline{[k]}_{\GG}$. In other words, images $\overline{[\O]}_\GG$ and $\overline{[\O]}_{\GK}$ are conjugate by $\overline{[\Phi]}$. Let $\eta\in \Aut(\GO)$ such that $\eta\overline{[\O]}_\GG\eta^{-1} = \O$, $\widetilde{\eta}$ an lifting of $\eta$. We replace $\Phi$ by $\widetilde{\eta}\circ\Phi$ so $\overline{[\O]}_\GG=\O$.
\subsection{Equi-Height Pairs}
Note $\Phi$ induces a map $\biso=\overline{[\Phi]}$ on the special fiber
$$
\map{\biso}{\GG_0}{\GG_0}.
$$
Identifying $T\GK\cong\OO K$ and $T\GG\cong\OO F^2$, by rational Tate modules of $\GK$ and $\GG$, $\Phi$ induces a map $\blin$ 
$$
\map{\blin}{F^2}K.
$$
In other words, $\Phi$ induces a map $\biso$ and a map $\blin$ up to right-$\GL_2(\OO F)$ and left-$\OO K^\times$ action and
$$[\OO K:\tau(\OO F^2)] = \deg(\biso).$$
We call such a $(\biso,\blin)$ an equi-height pair. Conversely, any equi-height pair $(\biso,\blin)$ determines $\Phi$ up to *-isomorphisms. Here by *-isomorphism we mean an isomorphism that reduces to the identity map in the special fiber. An automorphism is a *-isomorphism only when it is the identity map. The data $(\biso,\blin)$ determines a (*)-isomorphic class of quasi-canonical liftings.

\subsection{The action of Galois group}
The action of $\Gal(\overline{K}/\breve{K})$ on $T\GG_K$ gives a map $$\Gal(\overline{K}/\breve{K})\rightarrow \OO K^\times.$$ By Lubin and Tate, this is a surjective. Note that $\O^\times$ preserves the submodule $T\GG$. Let $\Gamma\subset\Gal(\overline{K}/\breve{K})$ be the preimage of $\O^\times$, and let $W$ be the fixed field of $\Gamma$. Since $\Gamma$ preserves $T\GG$, $\GG$ is a formal $\O$-module over $W$. We fix the identification
\[equation]{\lb{Galois}
\Gal(\W/\breve{K})\;\;\cong\;\;\Gal(\overline{K}/\breve{K})/\Gamma\;\;\cong\;\;\OO K^\times/\O^\times.
}
 For $k\in \OO K^\times$, Let $(k)$ be the corresponding element in $\Gal(\W/\breve{K})$, write its action as $t\mapsto t^{(k)}$.
\ss{What is a deformation, what is the parameter for the deformation}
\subsection{The universal deformation of $\GG_0$}
A formal module is called a deformation of $\GO$ if it reduces to $\GO$. For example, $\GG$ and $\GK$ are deformations of $\GO$. By Lubin and Tate in \cite{lubin1966formal}, there is an universal deformation formal $\OO F$-module $\GU$ over $\CO F[[U]]$ for $\GG_0$. For any complete local $\CO F$-algebra $A$ with residue field $\CFF q$, assigning the varible $U$ in coefficients of $\GU$ to a topolygically nilpotent element $t\in A$ defines a deformation of $\GG_0$ over $A$. 

\[defi]{[U-section] We call the above $t$ the U-section of the deformation.
}

We denote this deformation as $\GG^t$. Suppose $\map\beta{\CO F[[U]]}{A}$ is the map sends $U$ to $t$. We see $\GG^t=\beta^*\GU$. The formal module $\GU$ is universal in the sense that every formal $\OO F$-module deformation of $\GG_0$ over $A$ is *-isomorphic to some $\GG^t$ for some $t\in A^\circ$. Here the set $A^\circ$ is the set of topologically nilpotnetn elements in $A$.

Therefore, Two height 2 formal $\OO F$-modules over $A$ are *-isomorphic if and only if their U-sections are the same. Since the choice of $U$ is not canonical, we need to fix the universal formal module $\GG^{univ}$. The set $A^\circ$ classifies all *-isomorphic classes of deformations of $\GG_0$ over $A$. 

\subsection{Automorphism Liftings of $\GG_0$}

Since the endomorphism ring of $\GO$ is the maximal order of a division algebra. For any $a\in\End(\GG_0)$, either $a$ or $1-a$ must be an automorphism. The problem of lifting endomorphisms reduces to lifting automorphisms of $\GO$. 

Let $\ACj\in\Aut(\GO)$, $[\ACj]_{\GO}$ represents a power series over $\CFF q$. Lifting $[\ACj]_{\GO}$ arbitrarily to a power series $\psi$ with coefficients in $\CO F[[U]]$. Since $\psi^\prime(X)\neq 0$, it is invertible. By putting $\GG^U=\GG^{univ}$ and
$$
X[+]_{\GG^\prime}Y=\psi(\psi^{-1}(X)[+]_{\GG^U}\psi^{-1}(Y))\qquad [l]_{\GG^\prime}(X)=\psi([l]_{\GG^U}(\psi^{-1}(X));\qquad l\in\OO F^\times
$$
we produced another formal $\OO K$-module $\GG^\prime$ with the isomorphism
$$
\map{\psi}{\GU}{\GG^\prime}.
$$
By an abuse of notation we denote the U-section of $\GG^\prime$ by $\ACj(U)$. There is a *-isomorphism 
$$
\map{\Phi}{\GG^\prime}{\GG^{\ACj(U)}}
$$
Let $\psi_{\ACj}=\Phi\circ\psi$ we obtained an isomorphism
\begin{equation}\lb{doubledating}
\map{\psi_\ACj}{\GU}{\GG^{\ACj(U)}}.
\end{equation}
The uniqueness of $\psi_{\ACj}$ and $\GG^{\ACj(U)}$ is followed by the universal property. The isomorphism \eqref{doubledating} is an universal lifting of $\ACj\in\Aut(\GO)$.

From now on, we view $\ACj(U)$ as a power series over $\CO F$ with the varible $U$. For any $a\in A^\circ$, $\ACj(a)$ represents applying this power series to $a$. This defines an action of $\Aut(\GO)$ on $A^\circ$. Please note that the power series $\ACj(U)$ we used in our discussion is different from the power series $[\ACj]_{\GO}$ over $\CFF q$ defining automorphisms of $\GO$.
From the universal lifting \eqref{doubledating}, we know for any $t,s\in A^\circ$, $\ACj$ can be lifted as an isomorphism $\map{\psi_{\ACj}}{\GG^t}{\GG^s}$ if and only if $s=\ACj(t)$. In particular, $\ACj$ can be lifted as an automorphism of $\GG$ if and only if $\ACj$ fixes its U-section. 

For any deformation $\GG$ of $\GO$ over a discrete valuation ring $A$ with uniformizer $\varpi$ and an arbitrary automorphism $\ACj\in\Aut(\GO)$. Let $t$ be the U-section for $\GG$. the maximal $n$ such that $\ACj$ can be lifted as an automorphism of $\GG$ over $A/\varpi^n$ is the maximal $n$ that makes \[equation]{\lb{yinyue}\ACj(t)\equiv t\text{ mod }\varpi^n.}
Therefore $n$ equals the valuation of $t-\ACj(t)$.

\ss{Action of the maximal order of division algebra and action of Galois group on parameters.}

\[prop]{\lb{Galo}Let $t$ be the U-section of $\GG$, if $k\in\O\left[\frac1\pi\right]\cap\OO D^\times$, then
$$
k(t) = t^{(k)}
$$
}
\[proof]{By \eqref{youyin}, as submodules of $\OO D$, we have $\O = \biso[\O]_{\GK}\biso^{-1}$. We need to find an isomorphism $\map\alpha\GG{\GG^{(k)}}$ that reduces to $\biso\overline{[k]}_{\GK}\biso^{-1}$ in $\GO$. Conjugating the isogeny in \eqref{fangong} we have
$$
\map{\Phi^{(k)}}{\GK}{\GG^{(k)}}.
$$
The map $\Phi^{(k)}$ and $\Phi$ embeds $T\GG$ and $T\GG^{(k)}$ into $T\GK$. By definition of $(k)$, we have 
$$
T\GG^{(k)} = k\cdot T\GG
$$
as submodules of $T\GK$. Note the endomorphism $[k]_{\GK}$ multiplies $k^{-1}$ on $T\GK$. The isogenies
$$
\map{\Phi}{\GK}{\GG};\qquad
\map{\Phi^{(k)}\circ [k]_{\GK}}{\GK}{\GG^{(k)}}
$$  
will send $T\GG$ and $T\GG^{(k)}$ to the same subset of $T\GK$. Therefore there exists an isomorphism $\map{\alpha}\GG{\GG^{(k)}}$ such that
$$
\alpha\circ\Phi=\Phi^{(k)}\circ[k]_{\GK}.
$$
Since $\Phi$ and $\Phi^{(k)}$ reduces to $\biso$ on the special fiber, $\alpha$ reduces to $\biso\overline{[k]}_{\GK}\biso^{-1}$.}

\subsection{Intersection numbers in Lubin-Tate space}
We put
$$
\DEF=\Spf(\CO F[[U]])
$$
and call it the Lubin-Tate deformation space of $\GG_0$.

The formal $\O$-module $\GG$ we defined in \eqref{fangong} corresponds to an assignment of $U$ to its U-section in $\W$, which give rise to a formal curve $\y$ in $\DEF$:
\begin{equation}\lb{xiangxiang}
\map{\y}{\Spf \OO W}{\DEF}
\end{equation}
In Lemma \ref{emptyset} we will prove this is a colosed embedding. The Lubin-Tate space $\DEF$ admits an action by $\Aut(\GO)$. By an abuse of notation, for any $\ACj\in\Aut(\GO)$, we denote $\map{\ACj}{\DEF}{\DEF}$ the isomorphism defined by $U\mapsto \ACj(U)$. Suppose $t$ is the U-section of $\GG$. Then $\ACj(\y)$ corresponds to $\GG^{\ACj(t)}$. We are interested in the valuation of $t-\ACj(t)$, which is the maximal $n$ that $\ACj$ can be lifted as an automorphism over $W/\varpi^{n}$. Put $$
\begin{array}{lll}
\DEF_W&=&\DEF\times_{\OO F}\Spf \OO W;\\
\DEF_K&=&\DEF\times_{\OO F}\Spf \CO K.\\
\end{array}
$$Let $\map{\x}{\Spf\OO W}{\DEF_W}$ be the graph of $\y$ in \eqref{xiangxiang}, $\map{\z}{\Spf\OO W}{\DEF_K}$ the map by composing $x$ with the projection to $\DEF_K$. Let $\a$ be the map by twisting $\z$ with the non-trivial Galois conjugate on the second factor.
Then $\x$ is a closed embedding. With an abuse of notation, denote by $\OO W$ the structural sheaf of $\Spf \OO W$. For $p=\x,\y,\z,\a$, define coherent sheaves on $\DEF_?=\DEF_W$,$\DEF_K$,$\DEF$ by
$$
\OO p=p_*\OO W.
$$
Let $y_1,y_2$ be different maps as \ref{xiangxiang}, we define
$$
\langle y_1,y_2 \rangle = \length_{\CO K}\left(\OO{y_1}\otimes_M\OO{y_2}\right).
$$
Furthermore, note that $\DEF_?$ admits an action by $\OO D^\times$, for any $\ACj\in\OO D^\times$ We define

\[equation]{\lb{tanamy}
\left\{
\begin{array}{llll}
\va_\x(\ACj)&=&\length_{\CO F}(\OO \x\otimes_{\DEF_W} \OO{\ACj \x}). \\
\va_\y(\ACj)&=&\length_{\CO F}(\OO \y\otimes_{\DEF} \OO{\ACj \y}). \\
\va_\z(\ACj)&=&\length_{\CO F}(\OO \z\otimes_{\DEF_K} \OO{\ACj \z}). \\
\va_\a(\ACj)&=&\length_{\CO F}(\OO \a\otimes_{\DEF_K} \OO{\ACj \z}). \\
\end{array}\right.
}
Let $\GG_{n}=\GG\otimes \OO W/\varpi^{n+1}$.  By interpretations in \eqref{yinyue}, we have
\begin{equation}\lb{xiangjiao}
\ACj\in\Aut(\GG_n)\Longleftrightarrow n<\va_x(\ACj).
\end{equation}
The proof of the theorem amounts to calculate $\va_\x(\ACj)$. 
\section{Computation of Intersection Numbers}

%% file: Lifting/proofs.tex
We will calculate $\va_\y(\ACj)$ and $\langle y_1,y_2\rangle$ in this section. We recommend the reader to skip but only remember Theorem \ref{shenqi} and Lemma \ref{landuo}. We use $|\bullet|_D$ and $|\bullet|_F$ to denote the absolute values for $D$ and $F$ such that $|\pi|_F=q^{-1}$ and $|\pi|_D=q^{-2}$.
\subsection{Results}
Let $\dd x$ be the normalized Haar measure for $\OO F$. For any $\ACj\in\End(\GO)$, let
\begin{equation}\lb{fyf}
\phi(\ACj)=\int_{\pi\OO F}|x-\ACj|_D^{-1}\dd x+1.
\end{equation}
Let $y_1,y_2$ be maps as \eqref{xiangxiang} corresponding to quasi-canonical liftings of formal $\OO 1,\OO2$-modules $\GG_1$, $\GG_2$ over $\OO{W_1}$ and $\OO{W_2}$, with U-sections $t_1$,$t_2$. Let $\mu_1,\mu_2$ be generators of $\OO 1$,$\OO 2$ with smallest absolute value. Choose uniformizers $\varpi_1,\varpi_2$ in $\OO{W_1}$ and $\OO{W_2}$. Let $\pi_D$ be an uniformizer of $\OO D$.
\[thm]{\lb{shenqi}
We have
\[equation]{\lb{qiguaizixi}
\langle y_1,y_2\rangle = \[cases]{\phi(\mu_1)&\text{ if } 1>|\mu_1|_D>|\mu_2|_D;\\1&\text{ if } 1=|\mu_1|_D>|\mu_2|_D.}
}
Let $\dd k$ be the normalized Haar measure of $\OO K^\times$. For any $\ACj\in\End(\GO)$, we have
\[equation]{\lb{xishoujian}
\va_y(\ACj)=\[cases]{\infty&\text{ if }\ACj\in \OO D^\times\cap(D^+\cup D^-);\\
[\OO K^\times:\O^\times]\left(1+\int_{\OO K^\times}|k-\ACj|_D^{-1}\dd k\right)&\text{ if }K/F\text{ ramifeid, }\ACj\in\OO K^\times+\pi_D\OO D;\\
[\OO K^\times:\O^\times]\int_{\OO K^\times}|k-\ACj|_D^{-1}\dd k&\text{ if }K/F\text{ unramifeid and }\OO K\neq\O;\\
\frac{r+1}2&\text{ if }K/F\text{ unramifeid and }\OO K=\O.
}}
Here $r=\va_D(\ACj_-)$, where $\ACj_-=(\overline\mu-\mu)^{-1}(\ACj\mu-\mu\ACj)$. The group $\OO D^\times\cap(D^+\cup D^-)$ is the normalizer of $\OO K^\times$ in $\OO D^\times$, where $D^+=K\subset D$ and $D^-=\{\gamma\in D:\gamma\mu=\overline\mu\gamma\}$.
}
The value $\langle y_1,y_2 \rangle$ leads the following observation
\[lem]{\lb{emptyset}
The map $\y$ in \eqref{xiangxiang} is a closed embedding.
}
\[proof]{We need to show the induced map 
\[equation]{\lb{tina}
\CO F[[U]]\longrightarrow\OO W = \OO y
}
is onto. If $|\mu|_D=0$, then $W=\breve F$ implies this map is surjective. When $|\mu|_D<1$, Let $\y_0$ be the map in \eqref{xiangxiang} corresponding to $\GK$, and $t_0$ its U-section. by Theorem \ref{shenqi}, we have
$$
\length(\OO{\y_0}\otimes_\DEF\OO y)=1.
$$ 
Since $\CO F[[U]]\longrightarrow\OO{y_0}$ is a surjective. Let $t$ be the U-section of $\GG$, the above expression implies
$$
\length\left(\OO y/(t-t_0)\right)=1.
$$
Since $W/\breve F$ is a ramified extension, so $t_0$ is not a uniformizer of $W$, therefore $t$ is a uniformizer, the image of $U$ is a generator of $\OO W$. The map \eqref{tina} is a surjective. We proved this Lemma. 
}
The value $\langle y_1,y_2 \rangle$ will also help us to determine $\va_\x(\ACj)$ in the following sense.
\[lem]{\lb{landuo}Let $\x_1$, $\x_2$ be graphs of $\y_1,\y_2$. If $\va_{\x_1}(\ACj)<\langle y_1,y_2 \rangle$, then
$$
\va_{\x_1}(\ACj)=\va_{\x_2}(\ACj).
$$
}
\[proof]{For $i=1,2$, since $y_i$ are closed embeddings. The natural maps
$$
\map{\iota_i}{\OO{W_i}}{\OO{W_1}\otimes_\DEF \OO{W_2}} 
$$
are surjective. This induces an isomorphism of coimages of $\iota_1$ and $\iota_2$
$$
\map{\alpha}{\OO{W_1}/\varpi_1^a}{\OO{W_2}/\varpi_2^a}
$$
for any $a\leq\length(\OO 1\otimes_\DEF\OO 2) = \langle y_1,y_2 \rangle$ such that $\alpha(\overline{t_1})=\overline{t_2}$. Since $\overline{t_1},\overline{t_2}$ are U-sections of $\GG_1\otimes\OO{W_1}/\varpi_1^a$ and $\GG_2\otimes\OO{W_2}/\varpi_2^a$,
\[equation]{\lb{xihuantina}
\End(\GG_1\otimes\OO{W_1}/\varpi_1^a) = \End(\GG_2\otimes\OO{W_2}/\varpi_2^a)
}
as subrings of $\End(\GO)$. Let $n=\va_{\x_1}(\ACj)$, we have
$$
\ACj\in \End(\GG_1\otimes\OO 1/\varpi_1^n)\smallsetminus\End(\GG_1\otimes\OO 1/\varpi_1^{n+1}).
$$
Since $n+1\leq\langle y_1,y_2 \rangle$, therefore by \eqref{xihuantina}
$$
\ACj\in \End(\GG_2\otimes\OO 2/\varpi_2^n)\smallsetminus\End(\GG_2\otimes\OO 2/\varpi_2^{n+1}).
$$
This implise $\va_{\x_2}(\ACj)=n=\va_{\x_1}(\ACj)$.}
\subsection{The intersection formula}We prove Theorem \ref{shenqi} by computation. Before further elaboration, we use $\dd g$ and $\dd k^\times$ for Haar measures on groups $\GL_2(\OO F)$ and $\OO K$ normalized by the subset $\GL_2(\OO F)$ and $\OO K^\times$ respectively. For any subset $A,B$ of those groups, by $[A:B]$ we mean
$$
[A:B] = \Vol(A)/\Vol(B).
$$
By an abuse of notation, we reserve $[W:\breve{K}]$ to denote the degree of the field extension $W/\breve{K}$. 

Let $[\OO y]$ be the class of $\OO y$ in the $\QQ$-coefficient K-group of $\DEF$. To offset the influence of $W$, we normalize it by
\[equation]{\lb{yaofabiao}
\delta[\biso,\blin]=\frac1{[W:\breve K]}\left[\OO y\right].
}
Let $(\biso_1,\blin_1)$ and $(\biso_2,\blin_2)$ be equi-height pairs for $\GG_1,\GG_2$. The intersection formula in \cite{qirui2017love} gives
\begin{equation}\label{gongshizi}\comments{gongshizi}
\chi(\delta[\biso_1,\blin_1]\otimes_{\DEF}\delta[\biso_2,\blin_2])=\frac{\zeta_{K_1}(1)\cdot\zeta_{K_2}(1)}{|\Delta_{K_1/F}|_F\cdot \zeta_F(1)\cdot\zeta_F(2)}\int_{\GL_2(\OO F)}|R(g)|_D^{-1}\dd g.
\end{equation}
Here
\begin{itemize}
\item $\zeta_L(s)=(1-q_L^{-s})^{-1}$, where $q_L$ is the residue cardinality of $\OO L$ for $L=F,K_1,K_2$.
\item $\Delta_{K_1/F}$ is the discriminant of $K_1/F$.
\item $\chi(\mathcal{F},\mathcal{G})$ represents the following number for coherent sheaves $\mathcal{F}$,$\mathcal{G}$ over $\DEF$:
$$
\chi(\mathcal{F},\mathcal{G})=\sum_{i=0}^{\infty}\length_{\CO F}\(\mathrm{Tor}^{(i)}_{\DEF}(\mathcal{F},\mathcal{G})\).
$$
\item The function $R(g)\in\OO D$ depends on $\biso_1,\blin_1,\biso_2,\blin_2$ is the following expression
$$
R(g)=\begin{pmatrix}0&1\end{pmatrix}\begin{pmatrix}1&1\\\mu_1&\overline{\mu_1}\end{pmatrix}^{-1}\varphi_1^{-1}g\varphi_2\mm11{\mu_2}{{\overline\mu}_2}\cc10,
$$
\end{itemize}
We remark here that $\DEF$ is regular of dimension 2 and $\delta[\biso,\blin]$ is of dimension 1. Its higher Tor groups vanish (See Lemma 4.3 of \cite{qirui2017love} or \cite[\href{http://stacks.math.columbia.edu/tag/0B01}{Tag 0B01}]{stacks-project}). So we have
\[equation]{\lb{zhaonvpengyou}
\chi(\delta[\biso_1,\blin_1],\delta[\biso_2,\blin_2])=\length_{\CO F}\(\delta[\biso_1,\blin_1]\otimes_{\DEF}\delta[\biso_2,\blin_2]\).
}
Notice that for $i=1,2$, $$[W_i:\breve{K}]\zeta_{K_i}(1)=[\OO{K_i}^\times: \O_i^\times][\OO{K_i}:\OO{K_i}^\times] = [\OO{K_i}:\O_i^\times].$$ Furthermore, write
$$
\ep_F=\zeta_F(1)^{-1}\zeta_F(2)^{-1}=(1-q^{-1})(1-q^{-2}).
$$
By \eqref{gongshizi},\eqref{zhaonvpengyou}, we have
\begin{equation}\lb{halu}
\langle y_1,y_2\rangle = \ep_F[\OO{K_1}:\O_1^\times][\OO{K_2}:\O_2^\times]|\Delta_{K_1/F}|_F^{-1}\int_{\GL_2(\OO F)}|R(g)|_D^{-1}\dd g.
\end{equation}
For any $g\in\GL_2(\OO F)$ we will use $g_{ij}$ to denote $i,j$-th entry of $g$, we define the following subset
\[equation]{\lb{kunnan}
\Gamma(a)=\{g\in\GL_2(\OO F)| g_{21}\in a\OO F\};
}
\[equation]{\lb{sjy}
\Omega(a)=\{g\in\GL_2(\OO F)| g_{21}\in a\OO F^\times\}.
}
Now we compute $\langle y_1,y_2\rangle$ when $|\mu_1|_D>|\mu_2|_D$ by the formula \eqref{halu}.

\subsection{Computation for $\langle y_1,y_2\rangle$}
Let $\ACj_0=\biso_1^{-1}\biso_2$. Write $\mu_2^{(\ACj_0)}:=\ACj_0\mu_2\ACj_0^{-1}$. We have
\begin{equation}\label{aru}\comments{aru}
R(g)=|(\mu_1-\overline{\mu}_1)^2|_F\cdot|g_{21}+g_{11}\overline{\mu_1}+g_{22}\mu_2^{(\ACj_0)}+g_{12}\overline{\mu_1}\mu_2^{(\ACj_0)}|_D^{-1}|\ACj_0|_D^{-1}.
\end{equation}
Since either $g_{21}$ or $g_{11}$ is a unit, together with $|\mu_1|_D>|\mu_2|_D$, we have
$$
|g_{21}+g_{11}\overline{\mu}_1|_D\geq |\mu_1|_D>|g_{22}\mu_2^{(\ACj_0)}+g_{12}\overline{\mu}_1\mu_2^{(\ACj_0)}|_D.
$$
Therefore 
\begin{equation}\label{budeng}
R(g)=|\ACj_0|_D^{-1}|\mu_1-\overline{\mu}_1|_D\cdot|g_{21}+g_{11}\overline{\mu}_1|_D^{-1}\end{equation}
We can write
$$
\int_{\GL_2(\OO F)}R(g)=|\ACj_0|_D^{-1}|\mu_1-\overline{\mu}_1|_D\left(\int_{\Omega(1)} |g_{21}+g_{11}\overline\mu_1|_D^{-1}\dd g + \int_{\Gamma(\pi)}|g_{21}+g_{11}\overline\mu_1|_D^{-1}\dd g\right)
$$
Note that $|a_{21}+a_{11}\overline\mu_1|_D=1$ over $\Omega(1)$, and $\Vol(\Omega(1))=(1+q^{-1})^{-1}$. Let $\dd g_{11}^\times,\dd g_{21},\dd g_{22}^\times$ be Haar measures of $\OO F^\times$, $\OO F$, $\OO F^\times$ respectively, then the measure $q\dd g_{11}^\times\dd g_{21}\dd g_{22}^\times$ is normalized Haar measure for $\Gamma(\pi)$. So we have
$$
\dd g = q\Vol\left(\Gamma(\pi)\right)\dd g_{11}^\times\dd g_{21}\dd g_{22}^\times = \frac1{1+q^{-1}}\dd g_{11}^\times\dd g_{21}\dd g_{22}^\times.
$$
Furthermore by the replacement $a=g_{11}\mapsto g_{21}g_{11}$ for the integrand over $\Gamma(\pi)$, we write
$$
\frac1{|\ACj_0|_D^{-1}|\mu_1-\overline{\mu}_1|_D}\int_{\GL_2(\OO F)}R(g)=\frac1{1+q^{-1}} + \frac1{1+q^{-1}}\int_{\pi\OO F}|a-\overline\mu_1|_D^{-1}\dd a=\frac{\phi(\mu_1)}{1+q^{-1}}.
$$
Therefore we have
\[equation]{
\langle y_1,y_2\rangle =(1-q^{-1})(1-q^{-2})[\OO K:\OO 1^\times][\OO K:\OO 2^\times]|\ACj_0|_D^{-1}|\mu_1-\overline\mu_1|_D|\Delta_{K_1/F}|_F^{-1}\frac{\phi(\mu_1)}{1+q^{-1}}. 
}

Since $|\ACj_0|_D^{-1}=|\mu_2\mu_1^{-1}|_D$, $|1-\mu_1^{-1}\overline\mu_1|_D=|\Delta_{K_1/F}|_F$, $[\OO K:\OO 2^\times]|\mu_2|_D=(1-q^{-1})^{-1}$, we have
$$
\langle y_1,y_2\rangle = \frac{[\OO K:\OO 1^\times](1-q^{-2})}{1+q^{-1}}\phi(\mu_1).
$$
If $|\mu_1|_D=1$, then $[\OO K:\OO 1^\times]=(1-q^{-2})^{-1}$,  $|a-\mu_1|_D=1$ for $a\in\pi\OO F$, therefore $\phi(\mu_1)=1+q^{-1}$,
$$
\langle y_1,y_2\rangle = \frac{\phi(\mu_1)}{1+q^{-1}}=1.
$$
Otherwise, if $|\mu_1|_D<1$, then $[\OO K:\OO 1^\times]=(1-q^{-1})^{-1}$, in this case
$$
\langle y_1,y_2\rangle = \phi(\mu_1).
$$ We finished the calculation of the formula.
\subsection{Computation for $\va_y(\ACj)$} Note that $\va_y(\ACj)=\langle y,\ACj y\rangle$, in this case, the formula \eqref{halu} specializes to the case $K_1=K_2$,$\mu_1=\mu_2=\mu$, $\biso_1^{-1}\biso_2=\ACj$, we have
\[equation]{\lb{kuaidian}
\va_y(\ACj)=\ep_F[\OO K:\O^\times]^2\int_{\GL_2(\OO F)}\left|\rr10\mm11\mu{\overline\mu}^{-1}\ACj g\mm11\mu{\overline\mu}\cc01\right|_D^{-1}\dd g.
}
We should first simplify the integrand 
\[equation]{\lb{addiction}
I_\ACj(g):=\left|\rr10\mm11{\overline\mu}\mu^{-1}\ACj g\mm11{\overline\mu}\mu\cc01\right|_D^{-1}.
}
This equals to
\[equation]{\lb{lianai}
I_\ACj(g)=|\overline\mu-\mu|_D|\mu\ACj(g_{11}+g_{12}\mu)-\ACj(g_{21}+g_{22}\mu)|_D^{-1}.
}
If $|g_{21}|_D>|\mu|_D$, we know $I_\ACj(g)=|\overline\mu-\mu|_D|g_{21}|_D^{-1}$ by strong triangle inequality. In other words,
\[equation]{\lb{jingshen}
I_\ACj(g)=|\overline\mu-\mu|_D|g_{21}|_D^{-1}\qquad\text{for }g\in\Omega(\pi^n), |\pi^n|_D>|\mu|_D.
}
Let $u$ be the maximal integer that $|\pi^{u-1}|_D>|\mu|_D$, then $u=w$ if $K/F$ is unramified or $u=w+1$ otherwise. We decompose the integral \eqref{kuaidian} for $\va_y(\ACj)$ into two parts.

\[equation]{\lb{parts}
\frac{\va_y(\ACj)}{\ep_F[\OO K:\O^\times]^2}=\left(\sum_{0\leq n<u}\int_{\Omega(\pi^n)}|\overline{\mu}-\mu|_D|\pi^n|_D^{-1}\dd g\right)+\int_{\Gamma(\pi^u)}I(g)\dd g.
} 
Denote former and later part as $P_s$ and $P_d$, in other words,
\[equation]{\lb{sl}
P_s=\sum_{0\leq n<u}\int_{\Omega(\pi^n)}|\overline{\mu}-\mu|_D|\pi^n|_D^{-1}\dd g,
}
\[equation]{\lb{dp}
P_d=\int_{\Gamma(\pi^u)}I(g)\dd g.
}
\subsubsection{Computation of $P_s$}
The group $\GL_2(\OO F)$ acts transitively on $\PP^1(\OO F/\pi^n)$, with $\Gamma(\pi^n)$ the stablizer of the point representing the submodule $\OO F/\pi^n\oplus\{0\}\subset\left(\OO F/\pi^n\right)^2$. So when $n\geq 1$ the volumn of $\Gamma(\pi^n)$ is $(1+q^{-1})^{-1}q^{-n}$, which is the reciprocal of the cardinality of $\PP^1(\OO F/\pi^n)$. Note that $\Omega(\pi^n)=\Gamma(\pi^n)\smallsetminus\Gamma(\pi^{n+1})$, this implies
\[equation]{\lb{bz}
\int_{\Omega(\pi^n)}\dd g=\left\{\begin{array}{ll}(1+q^{-1})^{-1}&\text{ if }n=0\\\\(1+q^{-1})^{-1}(q^{-n}-q^{-n-1})&\text{ if }n>0\\\end{array}\right.
}
Since $|\pi|_D=q^{-2}$, we can use \eqref{bz} to calculate \eqref{sl}
\[equation]{\lb{slres}
P_s=|\mu-\overline{\mu}|_D\left((1-q^{-1})^{-1}+\sum_{1\leq n<u}(1-q^{-1})^{-1}(q^{-n}-q^{-n-1})q^{2n}\right) = |\mu-\overline{\mu}|_D(1+q^{-1})^{-1}q^{u-1}
}
\subsubsection{Computation of $P_d$}
Denote the subgroup
$$
\Gamma_0(a)=\left\{g\in\GL_2(\OO F):\mm{g_{11}}{g_{12}}{g_{21}}{g_{22}}=\mm10{g_{21}}{g_{22}}, g_{21}\in a\OO F\right\}
$$
For any $k\in\O$, let $k_H\in\GL_2(\OO F)$ be the element such that
$$
k_H\cc1\mu=\cc1\mu k.
$$
Let $H=\{k_H|k\in\O ^\times\}$, then $H$ is a subgroup of $\GL_2(\OO F)$.
\[lem]{
We have a group decomposition
$$
\Gamma(\pi^u)=\Gamma_0(\pi^u)H;\qquad \Gamma_0(\pi^u)\cap H=\left\{\mm1{}{}1\right\}
$$
}
\[proof]
{For any $g\in\Gamma(\pi^u)$, Let $k=g_{11}+g_{12}\mu$ then $k\in \O^\times$ and 
$$
gk_H^{-1}\cc1\mu=\cc1{\frac{g_{21}+g_{22}\mu}{g_{11}+g_{12}\mu}}=\cc1{c+d\mu},
$$
here we put $\frac{g_{21}+g_{22}\mu}{g_{11}+g_{12}\mu}=c+d\mu$. Since $1,\mu$ form an $\OO F$-basis of $\O$, we have $$gk_H^{-1}=\mm1{}cd.$$
Notice that $c=(g_{11}+g_{12}\mu)^{-1}(g_{11}+g_{12}\overline\mu)^{-1}(g_{21} g_{11}+g_{12} g_{22}\mu\overline\mu)\in\pi^u\OO F$ because both $g_{21},\mu\overline\mu\in\pi^u\OO F$. Therefore $gk_H^{-1}\in \Gamma_0(\pi^u)$.
To compute $\Gamma_0(\pi^u)\cap H$, note for every $h\in\Gamma_0(\pi^u)$,
$$
h\cc1*=\cc1*,
$$
where $*$ denotes arbitrary element. Since $h=k_H$ for some $k\in\O^\times$, then
$$
h\cc1\mu=\cc k{\mu k}
$$
Therefore $k=1$, we proved $\Gamma_0(\pi^u)\cap H$ is trivial.}

From now on, we write $g=lk_H$ for every $g\in \Gamma(\pi^u)$ where $l,k_H$ are elements of $\Gamma_0(\pi^u),H$ corresponding to the decomposition $\Gamma(\pi^u)=\Gamma_0(\pi^u)H$. We write the integrand \eqref{addiction} as
$$
I_\ACj(g)=|\overline\mu-\mu|_D\left|\rr{\mu}{-1}\ACj lk_H\cc1\mu\right|_D^{-1} = |\overline\mu-\mu|_D\left|\rr{\mu}{-1}\ACj l\cc1\mu k\right|_D^{-1}
$$

Let $\dd k_H$, $\dd l$ be normalized Haar measures on $H$ and $\Gamma_0(\pi^u)$, then $\dd l\dd k_H$ is normalized Haar measure for $\Gamma(\pi^u)$. Therefore
$$
\dd g = \Vol(\Gamma(\pi^u))\dd l\dd k_H=(1+q^{-1})^{-1}q^{-u}\dd l\dd k_H.
$$ 
We can simplify $P_d$ as
$$
P_d=(1+q^{-1})^{-1}q^{-u}\int_{\Gamma_0(\pi^u)}|\overline\mu-\mu|_D\left|\rr{\mu}{-1}\ACj l\cc1\mu\right|_D^{-1}\dd l\int_{\O^\times}|k|_D^{-1}\dd k.
$$

Since $|k|_D=1$ for all $k\in\O^\times$, we drop the last factor. By calculating the matrix, we have
$$
I_\ACj(l)=|\overline\mu-\mu|_D|\mu\ACj-\ACj\mu(l_{22}+l_{21}\mu^{-1})|_D^{-1}
$$Consider an inclusion map
$$
\maps{\iota}{\Gamma_0(\pi^u)}{\OO K^\times}l{l_{22}+l_{21}\mu^{-1}}
$$
The pull-back of Haar measures on $\OO K^\times$ are Haar measures on $\Gamma_0(\pi^u)$. The image of $\iota$ is $\OO F^\times\oplus \zeta\OO F$, which has relative volumn $1$ in $\OO K^\times$ if $K/F$ is ramified or otherwise $(1+q^{-1})^{-1}$. Now we denote $k=l_{22}+l_{21}\mu^{-1}$ and $\dd k^\times$ be the $\iota$-pull back of the Haar measure on $\OO K^\times$. Then we write $P_d$ as

\begin{equation}\lb{shenme}
\begin{split}
P_d=(1+q^{-1})^{-1}q^{-u}|\overline\mu-\mu|_D\int_{\OO K^\times}|\mu\ACj-\ACj\mu k|_D^{-1}\dd k^\times&\qquad\text{ if }K/F \text{ is ramified}.\\
P_d=|\overline\mu-\mu|_Dq^{-u}\int_{\OO F^\times\oplus\zeta\OO F}|\mu\ACj-\ACj\mu k|_D^{-1}\dd k^\times&\qquad \text{ if }K/F \text{ is unramified}. 
\end{split}
\end{equation}
To simplify the integrand $|\mu\ACj-\ACj\mu k|_D^{-1}$, we think $D$ as a left K-algebra. The right multiplication of $\mu$ decomposes $D$ as eigenspaces $D^+$ of eigenvalue $\mu$, and $D^-$ of eigenvalue $\overline\mu$. For $\ACj\in D$, write $$\ACj=\ACj_++\ACj_-$$ where $\ACj_+\in D^+$, $\ACj_-\in D^-$. So $\ACj_+\mu=\mu\ACj_+$, $\ACj_-\mu=\overline\mu\ACj_-$. Note $\ACj_-^2\in F$ since it commutes with $\mu,\ACj_-$. Denote $\overline{\ACj}=\overline{\ACj_+}-\ACj_-$ we have $\overline{\ACj}\ACj=\overline{\ACj_+}\ACj_+-\ACj_-^2\in F$.

\[lem]{\lb{xiebuwan}Suppose $K/F$ is ramified. We have $P_d=\infty$ if and only if $\ACj\in \OO D^\times\cap(D^+\cup D^-)$.}
\[proof]{The equation for $k\in K$
\[equation]{\lb{wangzhan}\mu\ACj-\ACj\mu k=0
}
have a solution if and only if $\ACj^{-1}\mu\ACj=\mu k\in K$. Since the minimal polynomial of $\mu$ and $\ACj^{-1}\mu\ACj$ are the same, either $\ACj^{-1}\mu\ACj=\mu$ or $\overline\mu$. This is equivalent to say $\ACj\in D^+\cup D^-$.

If \eqref{wangzhan} does not have a solution, the integrand of $P_d$ is continous on a compact subset $\OO F^\times\oplus\zeta\OO F$, therefore convergent. If \eqref{wangzhan} have a solution, $\mu^{-1}\ACj^{-1}\mu\ACj\in \OO K^\times$. Therefore the integral
$$
S_n = |\overline\mu-\mu|_D\int_{\ACj^{-1}\mu\ACj\mu^{-1}+\pi^n\OO K^\times}|\mu\ACj - \ACj\mu k|_D^{-1}\dd k^\times = |\overline\mu-\mu|_D|\ACj\mu|_D^{-1}.
$$
Since $P_d$ has positive integrand and $\mu^{-1}\ACj^{-1}\mu\ACj+\pi^n\OO K^\times$ are disjoint subsets of $\subset \OO K^\times$ We conclude $P_d>S_1+S_2+\cdots = \infty$.
}
\[cor]{\lb{huaxue}If $\sigma\in\OO D^\times\cap (D^+\cup D^-)$, then $\va_y(\sigma)=\infty$.}
\[proof]{By Lemma \ref{xiebuwan}, $P_d=\infty$, therefore $\va_y(\sigma)=\ep_F[\OO K:\O^\times]^2(P_s+P_d)=\infty$.}

\[lem]{\lb{jinzhang}If $|2|_D=1$, then for any $\ACj\in D$,
\[equation]{\lb{xueba}
|\ACj|_D=\max\{|\ACj_+|_D,|\ACj_-|_D\}
}
}
\[proof]{ By triangle inequality, we only need to show $
|\ACj|_D\geq\max\{|\ACj_+|_D,|\ACj_-|_D\}
$. Since
$$
(\mu-\overline\mu)\ACj(\mu-\overline\mu)^{-1}=\ACj_+-\ACj_-,
$$
we have $|\ACj_++\ACj_-|_D=|\ACj_+-\ACj_-|_D$, therefore
$$
|2\ACj_\pm|_D=|(\ACj_++\ACj_-)\pm(\ACj_+-\ACj_-)|_D\leq\max\{|\ACj_++\ACj_-|_D,|\ACj_+-\ACj_-|_D\}=|\ACj|_D.
$$
Therefore the Lemma follows.
}
Note $\ACj(\mu-\overline{\mu})\ACj_-=(\mu-\overline{\mu})\ACj_-\overline{\ACj}$, so $|\overline{\ACj}|_D=|\ACj|_D$. In the integrand of \eqref{shenme}, $|\ACj|_D=|\overline{\ACj}|_D=1$,
\[equation]{\lb{tiaopi}
|\mu\ACj-\ACj\mu k|_D^{-1} = |\overline{\ACj}(\mu\ACj-\ACj\mu k)|_D^{-1} = \left|\left(\mu(\overline{\ACj}_+\ACj_+-\overline\ACj\ACj k)-\ACj_-^2\overline\mu\right)+\ACj_-\ACj_+(\overline\mu-\mu)\right|_D^{-1}
}
We have $\mu(\overline{\ACj}_+\ACj_+-\overline\ACj\ACj k)-\ACj_-^2\overline\mu\in D^+$ and $\ACj_-\ACj_+(\overline\mu-\mu)\in D^-$. From now on, we assume $q$ is odd. This implies $|2|_D=1$ and $|\mu-\overline{\mu}|_D=|\mu|_D$. Now assume $\ACj\in\OO K^\times+\pi_D\OO D$, let $a\in\OO K^\times$ such that $\ACj-a\in\pi_D\OO D$. Therefore $|\ACj-a|_D<1$ and by Lemma \ref{jinzhang} this implies
$$
\max\left\{|\ACj_+-a|_D,|\ACj_-|_D\right\}<1.
$$  
Therefore $|\ACj_+-a|_D<1$, we have $|\ACj_+|_D=|a|_D=1$. The integrand simplified as
$$
|\mu\ACj-\ACj\mu k|_D^{-1}=\min\left\{|\mu(\overline\ACj_+\ACj_+-\overline\ACj\ACj k)-\ACj_-^2\overline\mu|_D^{-1},|\ACj_-\mu|_D^{-1}\right\}.
$$
Note that if $|\mu(\overline\ACj_+\ACj_+-\overline\ACj\ACj k)|_D^{-1}<|\ACj_-\mu|_D^{-1}$, we have
$$|\mu(\overline\ACj_+\ACj_+-\overline\ACj\ACj k)|_D^{-1}=|\mu(\overline\ACj_+\ACj_+-\overline\ACj\ACj k)-\ACj_-^2\overline\mu|_D^{-1}.$$
If $|\mu(\overline\ACj_+\ACj_+-\overline\ACj\ACj k)|_D^{-1}\geq|\ACj_-\mu|_D^{-1}$, since $|\ACj_-^2\overline\mu|_D\leq|\ACj_-\mu|_D$, we have
$$
|\mu(\overline\ACj_+\ACj_+-\overline\ACj\ACj k)-\ACj_-^2\overline\mu|_D^{-1}\geq |\ACj_-\mu|_D^{-1}.
$$
Therefore,
$$
\min\left\{|\mu(\overline\ACj_+\ACj_+-\overline\ACj\ACj k)-\ACj_-^2\overline\mu|_D^{-1},|\ACj_-\mu|_D^{-1}\right\}=
\min\left\{|\mu(\overline\ACj_+\ACj_+-\overline\ACj\ACj k)|_D^{-1},|\ACj_-\mu|_D^{-1}\right\}.
$$
So we can write
\[equation]{\lb{xiaozuozaina}
|\mu\ACj-\ACj\mu k|_D^{-1}=|\mu|_D^{-1}\min\left\{|\overline\ACj_+\ACj_+-\overline\ACj\ACj k|_D^{-1},|\ACj_-|_D^{-1}\right\}.
}
By replacing the variable $k\mapsto \overline\ACj_+(\overline\ACj\ACj)^{-1}k$ and using \eqref{xiaozuozaina} we write
\[equation]{\lb{ticket}
\int_{\OO K^\times}|\mu\ACj-\ACj\mu k|_D^{-1}\dd k=|\mu|_D^{-1}\int_{\OO K^\times}\min\left\{|\ACj_+-k|_D^{-1},|\ACj_-|_D^{-1}\right\}\dd k = |\mu|_D^{-1}\int_{\OO K^\times}|k-\ACj|_D^{-1}\dd k^\times.
}
Now we will use the formula of $P_s$ in \eqref{slres} and $P_d$ in \eqref{shenme} to calculate $\va_y(\ACj)$ in each cases.
\subsubsection{The ramified case}
If $K/F$ is ramified, We have \eqref{ticket} equals to $(1+q^{-1})|\overline\mu-\mu|_D^{-1}q^u P_d$,
$$
\frac{\va_y(\ACj)}{\ep_F[\OO K:\O^\times]^2|\Delta_{K/F}|_F^{-1}}=P_s+P_d=\frac{q^{-u}}{1+q^{-1}}+\frac{q^{-u}}{1+q^{-1}}\int_{\OO K^\times}|k-\ACj|_D^{-1}\dd k.
$$
Since $\ep_F=(1-q^{-1})^2(1+q^{-1})$, $[\OO K:\OO K^\times]=(1-q^{-1})^{-1}$, $[\OO K:\O^\times]|\Delta_{K/F}|_F^{-1}=(1-q^{-1})^{-1}q^u$,
\[equation]{\lb{ramified}
\va_y(\ACj)= [\OO K^\times:\O^\times]\left(1+\int_{\OO K^\times}|k-\ACj|_D^{-1}\dd k\right).
} 
\subsubsection{Unramified cases}If $K/F$ is unramified, we note that $\OO K^\times = \left(\OO F^\times\oplus\zeta\OO F\right)\coprod\left(\pi\OO F\oplus\zeta\OO F^\times\right)$. By equation \eqref{shenme}, for $|\overline\mu-\mu|_D^{-1}q^u P_d$, the equation \eqref{ticket} has overcounted the part
\[equation]{\lb{future}
|\mu|_D^{-1}\int_{\pi\OO F\oplus\zeta\OO F^\times}\min\left\{|\overline\ACj_+\ACj_+-\overline\ACj\ACj k|_D^{-1},|\ACj_-|_D^{-1}\right\}\dd k.
}
Since $\overline\ACj_+\ACj_+,\overline\ACj\ACj\in\OO F^\times$, we have $\overline\ACj_+\ACj_+-\overline\ACj\ACj k\in\OO F^\times\oplus\zeta\OO F^\times$. The integrand is 1. \eqref{future} equals to
$$
 |\mu|_D^{-1}\Vol(\pi\OO F\oplus\zeta\OO F^\times)=|\mu|_D^{-1}\frac{1}{1+q}.
$$

If $|\mu|_D<1$, this equals to $q^u P_s$ therefore we have
$$
\frac{\va_y(\ACj)}{\ep_F[\OO K:\O^\times]^2}=P_s+P_d=q^u\int_{\OO K^\times}|k-\ACj|_D^{-1}\dd k.
$$
Since $\ep_F=(1-q^{-1})^2(1+q^{-1})$, $[\OO K:\OO K^\times]=(1-q^{-2})^{-1}$, $[\OO K:\O^\times]=(1-q^{-1})^{-1}q^u$,
\[equation]{\lb{unramified}
\va_y(\ACj)= [\OO K^\times:\O^\times]\int_{\OO K^\times}|k-\ACj|_D^{-1}\dd k.
}
If $|\mu|_D=1$, by \eqref{shenme}, \eqref{xiaozuozaina} and change varible $k\mapsto \overline\ACj_+\ACj_+(\overline\ACj\ACj)^{-1}k$ we have directly
$$
\frac{\va_y(\ACj)}{\ep_F[\OO K:\OO K^\times]^2}=P_d=(1+q^{-1})\int_{\OO F^\times\oplus\zeta\OO F}\min\left\{|1-k|_D^{-1},|\ACj_-|_D^{-1}\right\}\dd k^\times.
$$
Note that$|1-k|_D^{-1} = q^{2n}$ if and only if $k\in 1+\pi^n\OO K^\times$. Let $a=0.5(\va_D(\ACj_-)+1)$We write $P_d$ into
$$
\left(\Vol\left(\OO F^\times\oplus\zeta\OO F\right)-\Vol\left(1+\pi\OO K\right))\right)+\sum_{n=1}^{a-1}q^{2n}\Vol\left(1+\pi^n\OO K^\times\right)+\Vol\left(1+\pi^{a}\OO K\right)q^{2a-1}
$$ 
By $\Vol(\OO F^\times\oplus\zeta\OO F)=(1+q^{-1})^{-1}$, $\Vol(1+\pi^n\OO K^\times)=q^{-2n}$, $\Vol(1+\pi^n\OO K)=(1-q^{-2})^{-1}q^{-2n}$, 
$$
P_d = (1+q^{-1})^{-1}\frac{\va_D(\ACj_-)+1}2.
$$
Since $\ep_F[\OO K:\OO K^\times]^2=(1+q^{-1})^{-1}$, we have
$$
\va_y(\ACj)=\frac{\va_D(\ACj_-)+1}2.
$$

%% file: Lifting/relation.tex
\section{Proof of main theorem}\lb{npy}
In this section, we will first define the notion of distance. Then prove the main theorem by classifying automorphisms of $\GG_0$ into two classes based on its distance to the subring $\O$. We call one class the shallow automorphism and another one deep automorphism. We prove the theorem for both of them with different methods by our formulae in Theorem \ref{shenqi}
\subsection{Distance and Projection}

\[defi]{
Let $\ACj\in D$, $S\subset D$ is a compact subset, The distance of $\ACj$ to $S$ is
$$
||\ACj||_S=\min_{x\in S}\{|\ACj-x|_D\}.
$$
Furthermore, we define the set of projection
$$
\Proj_S(\ACj)=\{a\in S: |\ACj-a|_D=||\ACj||_S\}.
$$
}

\[lem]{\lb{xiaozuo}
Let $\ACj\in D$, $S\subset D$ a compact subset, $\ACj^\prime\in\Proj_S(\ACj)$. \[enumerate]{
\item For any $a\in S$, $|\ACj+a|_D=\max\{|\ACj^\prime+a|_D,|\ACj-\ACj^\prime|_D\}$
\item If $\ACj^\prime\in T\subset S$, then $||\ACj||_T=||\ACj||_S$.
}
}
\[proof]{Since $|\ACj+a|_D=|(\ACj-\ACj^\prime)+(\ACj^\prime+a)|_D$, the equality holds if $|\ACj-\ACj^\prime|_D\neq |\ACj^\prime+a|_D$. Otherwise, if $|\ACj-\ACj^\prime|_D= |\ACj^\prime+a|_D$, 
$$
\max\{|\ACj^\prime+a|_D,|\ACj-\ACj^\prime|_D\} = |\ACj-\ACj^\prime|_D = ||\ACj||_S\leq |\ACj+a|_D.
$$
By triangle inequality, the last sign must be equal.  To prove (2), we note that both $||\ACj||_T$ and $||\ACj||_S$ equals to $|\ACj-\ACj^\prime|_D$. Therefore, $||\ACj||_T=||\ACj||_S$.
}

\subsection{Lifting of shallow automorphisms}
\ss{We set up\[itemize]{
\item The defintiion of shallow automorphism
\item The direct computation of shallow automorphism
}
}

We call $\ACj\in\OO D^\times$ a shallow automorphism to $\O$ or an $\O$-shallow automorphism if
$$
||\ACj||_{\OO F^\times}\geq |\pi^{-1}\mu|_D.
$$ 
Shallow automorphisms exists only when $|\mu|_D<1$, in this case, the U-section of $\GG$ is an uniformizer of $\OO W$ by Lemma \ref{emptyset}.
\[lem]{\lb{Anna}Let $\ACj^\prime\in \Proj_{\OO F^\times}(\ACj)$, $\ACj^{\prime\prime}=\ACj-\ACj^\prime$. If $\ACj$ is a shallow automorphism,$$\phi(\ACj^{\prime\prime})\leq\phi(\pi^{-1}\mu)$$
}
\[proof]{We compare the integrand for \eqref{fyf}. Note that for any $x\in\pi\OO F$, using Lemma \ref{xiaozuo} we have
$$|\ACj^{\prime\prime}-x|_D=|\ACj-\ACj^\prime-x|_D=\max\{|x|_D,||\ACj||_{\OO F^\times}\}.$$
And we have
$$|\pi^{-1}\mu-x|_D=\max\{|x|_D,|\pi^{-1}\mu|_D\}$$
The lemma follows since we defined $||\ACj||_{\OO F^\times}\geq |\pi^{-1}\mu|_D$.
}
\[lem]{\lb{fendou} If $\ACj\in\OO D$ and $\ACj\notin\OO D^\times$, we have $\phi(\pi\ACj)=q\phi(\ACj)$.}
\[proof]{Since $\ACj\notin\OO D^\times$ we have $|x-\pi\ACj|_D^{-1}=q^2$ for $x\in\pi\OO F^\times$, therefore
$$
\phi(\pi\ACj)=1+\int_{\pi\OO F}|x-\pi\ACj|_D^{-1}\dd x=1+\int_{\pi\OO F^\times}q^2\dd x+\int_{\pi\OO F}|\pi x-\pi\ACj|_D^{-1}\dd \pi x.
$$
Note that the subset $\pi\OO F^\times$ have volumn $q^{-1}(1-q^{-1})$, the above equation equals to
$$
1+(q-1)+|\pi|_D^{-1}|\pi|_F\int_{\pi\OO F}|x-\ACj|_D^{-1}\dd x=q+q\int_{\pi\OO F}|x-\ACj|_D^{-1}\dd x.
$$
This is exactly the value of $q\phi(\ACj)$.
}
From now on $\dd k^\times$ is always the Haar-measure normalized by $\OO K^\times$, $t$ is the U-section of $\GG$.
\[lem]{\lb{sunjunyi}Suppose $\ACj\in \O\left[\frac1\pi\right]\cap \OO D$, but $\ACj\notin \O$, then
$$
\va_x(\ACj) = [\OO K^\times:\O^\times]\int_{\O^\times}|k-\ACj|_D^{-1}\dd k^\times.
$$
}
\[proof]{By Proposition \ref{Galo}, we have $\ACj(t)=t^{(\ACj)}.$Therefore,
\[equation]{\lb{fengyifan}
\va_{W}(\ACj(t)-t)=\va_{W}\left(t^{(\ACj)}-t\right).
}
By Lemma \ref{emptyset}, $t$ is an uniformizer of $\OO W$, therefore the expression \eqref{fengyifan} is the lower numbering of the element $(\ACj)$. The value of of \eqref{fengyifan} does not depend on the choice of uniformizer. We construct another uniformizer as following. We embed $\OO K$ into $\OO D$ such that it has the image $\OO D\cap \O\left[\frac1\pi\right]$ and construct the corresponding canonical lifting $\GG^\prime$. Let $W^\prime/\breve K$ be the field extension by adding all $\pi^m$ torsions of $\GG^\prime$. Assume here $m$ is large enough so $W\subset W^\prime$. 
Let $X$ be an $\pi^m$-torsion of $\GG^\prime$ which is not killed by $\pi^{m-1}$. Note $X$ is an uniformizer of $\OO{W^\prime}$ and $$\Gal(W^\prime/\breve{K})\cong\OO K^\times/1+\pi^m\OO K.$$ By previous construction in \S\ref{yuanwang}, $W$ is the subfield fixed by $\O^\times$.
Therefore the Norm $$Y=\mathrm{N}_{W/\breve{K}}(X)=\prod_{\sigma\in\Gal(W^\prime/W)} X^{(\sigma)}=\prod_{k\in\O^\times/1+\pi^m\OO K}[k]_{\GG^\prime}(X)$$
is an uniformizer of $\OO{W}$.Put
\[equation]{\lb{liuxinran}
f(x) = \prod_{k\in\O^\times/1+\pi^m\OO K}[k]_{\GG^\prime}(x)-Y 
}
Since $f([k]_{\GG^\prime}X)=0$ for all $k\in \O^\times/1+\pi^m\OO F$, we have the following expansion
\[equation]{\lb{dazuo}
f(x) = u(x)\prod_{k\in\O^\times/1+\pi^m\OO K}\left(x[-]_{\GG^\prime}[k]_{\GG^\prime}X\right).
}
Let $x=0$ in \eqref{dazuo} and compare valuations, we see $u(0)$ is an unit. Plug $x=X^{(\ACj)}$ into \eqref{liuxinran},
$$
f\left(X^{(\ACj)}\right) = \prod_{k\in\O^\times/1+\pi^m\OO K}[k]_{\GG^\prime}(X^{(\ACj)})-Y = Y^{(\ACj)}-Y.
$$
On the other hand, using expansion \eqref{dazuo}, we have
$$
f\left(X^{(\ACj)}\right) = u(X^{(\ACj)})\prod_{k\in\O^\times/1+\pi^m\OO K}[\ACj-k]_{\GG^\prime}X.
$$
Since $\va_{W^\prime}([\ACj-k]_{\GG^\prime}X)=|\ACj-k|_D^{-1}$, $\va_{W^\prime}\left(u(X^{(\ACj)})\right)=0$, and $\va_{W^\prime}=[W^\prime:W]\va_{W}$, we have
$$
\va_{W}\left(Y^{(\ACj)}-Y\right) = \frac1{[W^\prime:W]}\sum_{k\in\O^\times/1+\pi^m\OO K}|\ACj-k|_D^{-1}.
$$
Note that $[W^\prime:W]=[\O^\times:1+\pi^m\OO K]$, we write this sum into
\[equation]{\lb{sure}
\va_{W}\left(Y^{(\ACj)}-Y\right) = \int_{\O^\times}|\ACj-x|_D^{-1}\dd x^\times,
}
here $\dd x^\times$ is normalized by $\O^\times$, therefore $\dd x^\times = [\OO K^\times:\O^\times]\dd k^\times$. We proved this Lemma.}
Now we generalize this formula for all $\O$-shallow automorphisms.
\[lem]{\lb{jiaqi} If $\ACj$ is $\O$-shallow automorphism, let $\ACj^\prime\in\Proj_{\OO F^\times}(\ACj)$ and $\ACj^{\prime\prime}=\ACj-\ACj^\prime$, then

$$[\OO K^\times:\O^\times]\int_{\O^\times}|\ACj-k|_D^{-1}\dd k^\times = \frac q{q-1}\phi(\ACj^{\prime\prime})-\frac2{q-1}$$

and this value is smaller than $\phi(\mu)$.
}
\[proof]{
Let $\dd x^\times$ be the Haar measure normalized by $\O^\times$. Then $\dd x^\times=[\OO K^\times:\O^\times]\dd k^\times$. Write $\O^\times = \OO F^\times\oplus \mu\OO F$ and decompose $x=a+b\mu$, then $\dd x^\times=\dd a^\times\dd b$, where $\dd a^\times$,$\dd b$ are Haar measures for $\OO F^\times$ and $\OO F$. Since $\ACj$ is $\O$-shallow, we have for any $a\in\OO F^\times$,
$$|\ACj-a|_D\geq ||\ACj||_{\OO F^\times}>|\mu|_D.$$
Using the triangle inequality, the equation \eqref{sure} can be written as
\[equation]{\lb{whoisnext}
\int_{\O^\times}|\ACj-a-b\mu|_D^{-1}\dd x^\times = \int_{\OO F}\int_{\OO F^\times}\left|\ACj-a\right|_D^{-1}\dd a^\times\dd b=\int_{\OO F^\times}\left|\ACj-a\right|_D^{-1}\dd a^\times
}

By Lemma \ref{xiaozuo}, 
\[equation]{
\int_{\OO F^\times}\left|\ACj-a\right|_D^{-1}\dd a^\times=\int_{\OO F^\times}\min\{|\ACj^{\prime\prime}|_D^{-1},|(\ACj^\prime-a)|_D^{-1}\}\dd a^\times.
}
Suppose $\OO F^{\times\times}=\OO F^\times\smallsetminus(\ACj^\prime+\pi\OO F)$. We found $|\ACj^\prime-a|_D=1$ on $\OO F^{\times\times}$, therefore \eqref{xiaozuo} equals to
$$
\Vol(\OO F^{\times\times})+\int_{\ACj^\prime+\pi\OO F}|\ACj^{\prime\prime}+\ACj^\prime-a|\dd a^\times.
$$
Note that $\Vol(\OO F^{\times\times})=\frac{q-2}{q-1}$, $\dd a^\times = \frac{q}{q-1}\dd a$, change variable $a\mapsto a+\ACj^\prime$,the integral simplifies to
$$
\frac{q-2}{q-1}+\frac q{q-1}\int_{\pi\OO F}|\ACj^{\prime\prime}-a|\dd a=\frac q{q-1}\phi(\ACj^{\prime\prime})-\frac2{q-1}.
$$
This value is less than $\phi(\mu)$ because by Lemma \ref{Anna} and Lemma \ref{fendou},
$$
\frac q{q-1}\phi(\ACj^{\prime\prime})-\frac2{q-1}\leq\frac q{q-1}\phi(\pi^{-1}\mu)-\frac2{q-1}<\frac {1}{q-1}\phi(\mu)\leq\phi(\mu).
$$The Lemma follows.}

\[thm]{\lb{syr}Let $\ACj$ be an $\O$-shallow automorphism, then
$$
\va_x(\ACj) = [\OO K^\times:\O^\times]\int_{\O^\times}|k-\ACj|_D^{-1}\dd k^\times.
$$ 
}

\[proof]{ Let $K_2=F[\ACj]$ and $\O_2=\OO F[\pi^m\ACj]$ for a large enough $m$ such taht $\O_2\subset \O+\mu\OO D$. Let $W_2$ be the abelian extension of $K$ corresponding to $\O_2^\times$, and $\OO {W_2}$ is its ring of integers. Let $\GG_2$ be a formal $\O_2$-module over $\OO {W_2}$ and let $t_2$ be its U-section, $\mu_2=\pi^m\ACj$. Let $y_2$ be the map as constructed in \eqref{xiangxiang} for $\GG_2$, $x_2$ the graph of $y_2$, take $\ACj^\prime\in\Proj_{\OO F^\times}(\ACj)$ and $\ACj^{\prime\prime}=\ACj-\ACj^\prime$. Since $\ACj\in \OO {K_2}^\times$, by Lemma \ref{sunjunyi} and Lemma \ref{jiaqi}, we have
$$
\va_{x_2}(\ACj)=\frac q{q-1}\phi(\ACj^{\prime\prime})-\frac2{q-1}.
$$
By Lemma \ref{jiaqi}, this value is less than $\phi(\mu)$. Since $|\mu_2|_D=|\pi^m\ACj|_D<|\mu|_D<1$, by Theorem \ref{shenqi}, we have $\langle y,y_2\rangle=\phi(\mu)$. Therefore, $\va_{x_2}(\ACj)<\langle y,y_2\rangle$, by Lemma \ref{landuo} and Lemma \ref{jiaqi}, we have
$$
\va_{x}=\va_{x_2}(\ACj)=\frac q{q-1}\phi(\ACj^{\prime\prime})-\frac2{q-1}=[\OO K^\times:\O^\times]\int_{\O^\times}|k-\ACj|_D^{-1}\dd k^\times.
$$
This Theorem follows. This proves the main Theorem for $\O$-shallow automorphisms.}

\subsection{Lifting of deep automorphisms}
We call $\ACj\in\OO D^\times$ an $\O$-deep automorphism if
$$
||\ACj||_{\OO F^\times}< |\pi^{-1}\mu|_D.
$$ 
\[lem]{\lb{canojiu} Suppose $\ACj_1,\ACj_2\in\OO D^\times$ and $||\ACj_1||_{\OO F^\times}>||\ACj_2||_{\OO F^\times}$, then $||\ACj_1\ACj_2||_{\OO F^\times}\geq||\ACj_1||_{\OO F^\times}$.}
\[proof]{Let $\ACj_1^\prime\in\Proj_{\OO F^\times}(\ACj_1)$, $\ACj_2^\prime\in\Proj_{\OO F^\times}(\ACj_2)$. For any $a\in\OO F^\times$, we can write
$$
\ACj_1\ACj_2-a=\ACj_1\ACj_2^\prime-a+\ACj_1(\ACj_2-\ACj_2^\prime).
$$
On one hand, since $\ACj_2^\prime\in\OO F^\times$,
$$
|\ACj_1\ACj_2^\prime-a|_D=|\ACj_1-\ACj_2^{\prime-1}a|_D\geq ||\ACj_1||_{\OO F^\times}.
$$
On the other hand, since $\ACj_1\in\OO D^\times$, and $||\ACj_2||_{\OO F^\times}< ||\ACj_1||_{\OO F^\times}$,
$$
|\ACj_1(\ACj_2-\ACj_2^\prime)|_D=|\ACj_2-\ACj_2^\prime|_D=||\ACj_2||_{\OO F^\times}<||\ACj_1||_{\OO F^\times}.
$$
So $|\ACj_1\ACj_2-a|_D\geq ||\ACj_1||_{\OO F^\times}$ for any $a\in\OO F^\times$.
}

\[cor]{\lb{jiandu}
Suppose $\ACj_1$ and $\ACj_2$ are $\O$-shallow and $\O$-deep automorphisms respectively. $\ACj_1\ACj_2$ is an $\O$-shallow automorphism.
}
\[proof]{ we have $||\ACj_2||_{\OO F^\times}<|\pi^{-1}\mu|_D\leq ||\ACj_1||_{\OO F^\times}$, therefore $||\ACj_1\ACj_2||_{\OO F^\times}\geq|\pi^{-1}\mu|_D$ by Lemma \ref{canojiu}.}

We will compute $\va_\x(\ACj)$ by $\va_\y(\ACj)$ for $\O$-deep automorphisms. Now let $t$ be the U-section of $\GG$, $v_W$ the normalized valuation in $W$. $g(T)\in\CO K[T]$ the minimal polynomial of $t$ over $\breve{K}$. 
\begin{lem}\lb{paobu}
We have
\begin{equation}\lb{xingxing}
v_{\z}(\ACj)=\sum_{k\in\OO K^\times/\O^\times}v_{\x}(k\ACj).
\end{equation}
\end{lem}
\begin{proof}
Since $\z$ is a closed embedding, we know $v_{\z}(\ACj)$ equals to the length of $$\frac{\CO K[[U]]}{\left(g(U)\right)+\left(g(\ACj(U))\right)}\cong \OO W/\(g(\ACj(t))\).$$ Therefore
$
v_{\z}(\ACj) = v_W\(g(\ACj(t))\)
$,
In contrast, $v_x(\ACj)$ equals to the length of 
$$\frac{\OO W[[U]]}{(U-t)+(\ACj(U)-t)}\cong \OO W/\ACj(t).$$ Therefore, 
$
v_{\x}(\ACj) = v_W(\ACj(t)-t)
$.
Note that $g(T)$ is the product of all $(T-t^{(k)})$ for $k$ run through $\OO K^\times/\O^\times$, therefore,
$$
\va_\z(\ACj)=v_W\(g(\ACj(t))\)=\sum_{k\in\OO K^\times/\O^\times} v_W(\ACj(t)-t^{(k)})=\sum_{k\in\OO K^\times/\O^\times}v_{\x}(k^{-1}\ACj),
$$
here the last equality is by $v_W(\ACj(t)-t^{(k)})= v_W(\ACj(t)^{(k^{-1})}-t)=v_W(k^{-1}\ACj(t)-t)=v_{\x}(k^{-1}\ACj)$.
\end{proof}
\subsection{Computation of $\va_{\z}$ if $K/F$ is unramified}
If $K/F$ is unramified, then $\va_{\y}=\va_{\z}$. Therefore 
\[equation]{\lb{ganqing}
\va_\z(\ACj)= [\OO K^\times:\O^\times]\int_{\OO K^\times}|k^{-1}\ACj-1|_D^{-1}\dd k^\times.
}

\subsection{Computation of $\va_{\z}$ if $K/F$ is ramified} If $K/F$ is ramified, Let $\sigma\in \OO D^\times$ be the element in the normalizer of $\O^\times$ but not centralizer of $\O^\times$. In other words, $\sigma\in \OO D^\times\cap D^-$.

\begin{lem}\lb{ximeng}
We have
\begin{equation}\lb{jieren}
\va_{\y}(\ACj)=\va_{\z}(\ACj)+\va_{\a}(\ACj).
\end{equation}
\end{lem}
\begin{proof}Since the minimal polynomial of $t$ over $\CO F$ is $g(T)\overline g(T)$, the sheaf of ideals for $\y$ is generated by $g(U)\overline g(U)$, for $\z$ is generated by $g(U)$, for $\a$ is generated by $\overline{g}(U)$. So  
\begin{equation}\lb{shala}
\va_{\y}(\ACj) = v_W\(g(\ACj(t))\overline{g}(\ACj(t))\)=v_W\(g(\ACj(t))\)+v_W\(\overline{g}(\ACj(t))\)=\va_{\z}(\ACj)+\va_{\a}(\ACj).
\end{equation}
This lemma follows.\end{proof}
\[lem]{\lb{shadoumei}We have
$\va_{\overline{z}}(\ACj)=\va_z(\ACj\sigma)$.}
\[proof]{
By Lemma \ref{xiebuwan}, $\va_y(\sigma)=\ep_F[\OO K:\O^\times]^2(P_s+P_d)=\infty$. We claim $\va_z(\sigma)<\infty$, if not, by formula \eqref{xingxing}, there exists $k\in\OO K^\times$ such that $\va_x(k\sigma)=\infty$. Therefore $k\sigma\in \Aut(\GG)$, so $\pi+k\sigma\in\Aut(\GG)$, which implies $\va_y(\pi+k\sigma)>\va_x(\pi+k\sigma)=\infty$. But $\pi+k\sigma\notin D_+\cup D_-$, contradiction by Lemma \ref{xiebuwan}. Therefore $\va_z(\sigma)<\infty$. Since $\infty=\va_\y(\sigma)=\va_\z(\sigma)+\va_\a(\sigma)$ but $\va_\z(\sigma)<\infty$, this implies $\va_\a(\sigma)=\infty$. This implies
$$
\va_W(\overline g(\sigma(t)))=\va_\a(\sigma)=\infty.
$$
So $\overline g(\sigma(t))=0$. Then the map $t\mapsto \sigma(t)$ lifts the non-trivial element in $\Gal(\breve K/\breve F)$ to $\Gal(W/\breve F)$,
$$
\va_\a(\ACj)=\va_W(\overline g(\ACj(t)))=\va_W(g(\ACj(\sigma(t))))=\va_\z(\ACj\sigma).
$$
This lemma follows}

\[cor]{\lb{zhaonvyou}
If $K/F$ ramified. Let $\epsilon\in\OO D^\times$ such that $||\epsilon||_{\OO F^\times}=1$, then
$$
\va_x(\epsilon) = 1.
$$
}
\[proof]{If $\O\neq \OO K$, $\epsilon$ is an $\O$-shallow automorphism, by formula in Theorem \ref{syr}, the integrand $|\epsilon-k|_D^{-1}=1$ for all $k\in\O^\times$. So $\va_x(\epsilon)=1$. If $\O=\OO K$, let $\mu_2=\pi\mu$, $\GG_2$ the quasi-canonical lifting as a formal $\OO F[\mu_2]$-module. Notice $$|\mu-a|_D^{-1}=|\mu|_D^{-1}=q$$ for all $a\in\pi\OO F$, we have $\phi(\mu)=1+q^{-1}q=2$. Since $||\epsilon||_{\OO F^\times}=1$, $\epsilon$ is a shallow automorphism of $\GG_2$. Let $y_2$ be corresponding map in \ref{xiangxiang} and $x_2$ its graph. Notice $|\epsilon-a-b\pi\mu|_D=|\epsilon-a|_D=1$ for any $x=a+b\pi\mu\in \OO 2^\times$, therefore
$$
\va_{x_2}(\epsilon) = [\OO{K_2}^\times:\O_2^\times]\int_{\O_2^\times}|\epsilon-k|_D^{-1}\dd k^\times = 1.
$$
Since $|\pi\mu|_D<|\mu|_D<1$ we have $\langle y,y_2 \rangle=\phi(\mu)=2$. By Lemma \ref{landuo},
$$
\va_x(\epsilon)=\va_{x_2}(\epsilon) = 1.
$$
The Corollary follows.}

Since $k\ACj\in\OO K^\times$ and $K/F$ is a ramified extension, we have $||k\ACj||_{\OO F^\times}<1$. Therefore $||\sigma||_{\OO F^\times}=1$ implies $||k\ACj\sigma||_{\OO F^\times} = 1$ by Lemma \ref{canojiu}. By Corollary \ref{zhaonvyou}, we have
$$
\va_x(k\ACj\sigma) = 1
$$ 
for all $k\in\OO K^\times/\O^\times$. By Equation \eqref{xingxing} and Lemma \ref{shadoumei}, we have
$$
\va_\a(\ACj)=\va_\z(\ACj\sigma)=\sum_{\OO K^\times/\O^\times}\va_x(k\ACj\sigma) = [\OO K^\times:\O^\times].
$$
Therefore, if $K/F$ is ramified, $\va_\z(\ACj)$ equals to
\[equation*]{\lb{lianai}
\va_\y(\ACj)-\va_\a(\ACj)=[\OO K^\times:\O^\times]\left(1+\int_{\O^\times}|\ACj-k|_D^{-1}\dd k\right)-[\OO K^\times:\O^\times] = [\OO K^\times:\O^\times]\int_{\O^\times}|\ACj-k|_D^{-1}\dd k^\times.
}
We found the expression of $\va_z(\ACj)$ in the ramified case is the same as unramifeid case \eqref{ganqing}.
\subsection{Proof of main Theorem}
Since $\va_z(\ACj)$ has the same expression in either cases, we compute $\va_\x$ uniformly, if $k\in\OO K^\times\smallsetminus\O$, then $k$ is an $\O$-shallow automorphism, so is $k\ACj$ by Corollary \ref{jiandu}. Therefore, $\va_x(\ACj)$ equals to
$$
\va_\z(\ACj)-\sum_{\substack{k\in\OO K^\times/\O^\times k\notin \O^\times}}\va_x(k^{-1}\ACj)=
[\OO K^\times:\O^\times]\left(\int_{\OO K^\times}|k^{-1}\ACj-1|_D^{-1}\dd k-\int_{\OO K^\times\smallsetminus\O^\times}|k^{-1}\ACj-1|_D^{-1}\dd k\right).
$$
This proves
$$
\va_x(\ACj)=[\OO K^\times:\O^\times]\int_{\O^\times}|\ACj-k|_D^{-1}\dd k^\times.
$$
Now we proved the Theorem for $\O$-deep automorphisms. The case for $\O$-shallow automorphisms was proved in Theorem \ref{syr}.